\documentclass[reqno]{amsart}

\usepackage{amsthm}
\usepackage{amssymb}
\usepackage{amsmath}

\newtheorem{theorem}{Theorem}[section]
\newtheorem*{theorem*}{Theorem}
\newtheorem{lemma}{Lemma}[section]

\newtheorem{proposition}{Proposition}[section]
\newtheorem*{proposition*}{Proposition}

\newtheorem{corollary}{Corollary}[section]
\theoremstyle{definition}
\newtheorem{definition}{Definition}[section]

\newtheorem*{remark*}{Remark}
\newtheorem*{remarks*}{Remarks}

\numberwithin{equation}{section}

\newcommand{\st}[1]{\ensuremath{^{\scriptstyle \textrm{#1}}}}

\newcommand{\CC}{{\mathbb C}}
\newcommand{\LL}{{\mathbb L}}
\newcommand{\RR}{{\mathbb R}}
\renewcommand{\SS}{{\mathbb S}}
\newcommand{\ZZ}{{\mathbb Z}}
\newcommand{\NN}{{\mathbb N}}

\newcommand{\fg}{{\mathfrak g}}

\newcommand{\ft}{{\mathfrak t}}

\newcommand{\D}{{\slash \! \! \! \partial }}
\newcommand{\spin}{\mbox{spin}}
\newcommand{\Res}{\mbox{Res}}
\newcommand{\E}{{\mathcal E}}

\renewcommand \a {\alpha}
\renewcommand \b {\beta}
\newcommand{\e}[1]{e^{2 \pi i#1}}

\begin{document}

\title[${\bf G}$-actions on graphs]
{$G$-actions on graphs}
\author[V. Guillemin]{V. Guillemin\footnotemark {*}}
\thanks{* Supported by NSF grant DMS 890771}
\address{Department of Mathematics, MIT, Cambridge, MA 02139}
\email{vwg@@math.mit.edu}
\author[C. Zara]{C. Zara\footnotemark {**}}
\thanks{** Supported by Alfred P. Sloan Doctoral Dissertation 
Fellowship grant DD 766}
\address{Department of Mathematics, MIT, Cambridge, MA 02139}
\email{czara@math.mit.edu}

\begin{abstract}
Let $G$ be an $n$-dimensional torus and $\tau$ a Hamiltonian
action of $G$ on a compact symplectic manifold, $M$.  If $M$ is 
pre-quantizable one can associate with $\tau$ a representation of 
$G$ on a virtual vector space, $Q (M)$, by $\spin^{\CC}$-quantization.  
If $M$ is a symplectic GKM manifold we will show that several 
well-known theorems about this ``quantum action'' of $G$:  for example, 
the convexity theorem, the Kostant multiplicity theorem and the 
``quantization commutes with reduction'' theorem for circle subgroups 
of $G$, are basically just theorems about $G$-actions on graphs.
\end{abstract}

\maketitle

\section{Introduction}
\label{sec:intro}

Let $\Gamma$ be a finite $d$-valent graph and let $G$ be a
$n$-dimensional torus.  In this paper we will be concerned with 
objects (rings, modules, $G$-representations $\ldots$)
associated to an ``action'' of $G$ on $\Gamma$.  To define what
we mean by this term, let $V_{\Gamma}=V$ be the vertices of $\Gamma$ 
and $E_{\Gamma}$ the oriented edges.  For each 
$e \in E_{\Gamma}$ let $i(e)$ and $t(e)$ be
the initial and terminal vertices of $e$, and let $\bar{e}$ be
the edge, $e$, with its orientation reversed.  (Thus
$i(e)=t(\bar{e})$ and $t(e)=i(\bar{e}$.)

\begin{definition}
Let $\varrho$ be a map which assigns to each oriented edge, $e$, of
$\Gamma$ a one dimensional representation, $\varrho_e$, with character
\begin{equation}
  \label{eq:1.1}
  \chi_e : G \to S^1 \, ,
\end{equation}
let $\tau$ be a map which assigns to each vertex, $p$, of
$\Gamma$ a $d$-dimensional representation, $\tau_p$, and let
$G_e$ be the kernel of (\ref{eq:1.1}).  $\varrho$ and $\tau$ define
an \emph{action} of $G$ on $\Gamma$ if they satisfy the axioms
(\ref{eq:1.2})--(\ref{eq:1.4}) below :
\begin{eqnarray}
  \label{eq:1.2}
  \tau_p & \simeq & \bigoplus_{i(e)=p} \varrho_e \\
\label{eq:1.3}
\varrho_{\bar{e}} & \simeq & \varrho^*_e\\
\label{eq:1.4}
\tau_{i(e)} |_{G_e} & \simeq  & \tau_{t(e)} |_{G_e} \, .
\end{eqnarray}
\end{definition}

\begin{remark*}
  For the connection between this graph-theoretic notion of
  ``$G$-action'' and the usual notion of $G$-action, see
  Section \ref{sec:gkm} (or, for more details, \cite[\S~3.1]{GZ2}).
\end{remark*}

Let $\ZZ^*_G$ be the weight lattice of $G$, and let $\alpha_e
\in \ZZ^*_G$ be the weight of the representation, $\varrho_e$, \emph{i.e.}
\begin{equation}
  \label{eq:1.5}
  \chi_e = e^{2\pi i \alpha_e} \, .
\end{equation}
By (\ref{eq:1.2}) and (\ref{eq:1.5}) both $\varrho$ and $\tau$ are
determined by the $\alpha_e$'s; so an action of $G$ on a graph,
$\Gamma$, can be thought of as a labeling of each edge, $e$, of
the graph by a weight, $\alpha_e$.  This labeling, however, will 
be forced by (\ref{eq:1.2})--(\ref{eq:1.4}) to satisfy certain
axioms.  For instance, by (\ref{eq:1.3})
\begin{equation}
  \label{eq:1.6}
  \alpha_{\bar{e}} = -\alpha_e \, .
\end{equation}
We will say that an action is a \emph{GKM action} if, for every
pair of edges with the same initial vertex, $p=i(e_1)=i(e_2)$,
\emph{either} $e_1=e_2$ or $\alpha_{e_1}$ and $\alpha_{e_2}$
\emph{are linearly independent}.  (For the geometric
interpretation of this property, see Section \ref{sec:gkm}).  All
the actions we consider below will be assumed to be GKM actions.

This paper is the fourth in a series of papers on the equivariant 
cohomology of graphs.  The first three papers in this series were 
concerned with the equivariant cohomology ring, $H_G(\Gamma)$.
In this paper we will be concerned with a slightly more
complicated object:  the equivariant ``$K$-cohomology'' ring of $
\Gamma$.  However, to motivate its definition, we will first
recall how $H_G(\Gamma)$ is defined:  Denote by $\fg$ and $\fg_e$ 
the Lie algebras of $G$ and $G_e$, and let $\SS(\fg^*)$ and
$\SS(\fg^*_e)$ be the symmetric algebras over the duals of $\fg$
and $\fg_e$.  From the inclusion of $\fg_e$ into $\fg$ one gets a 
restriction map
\begin{equation}
  \label{eq:1.7}
  r_e: \SS(\fg^*) \to \SS(\fg^*_e) \, .
\end{equation}

\begin{definition}
$H_G (\Gamma)$ is the set of all functions, $f:V_{\Gamma} \to \SS(\fg^*)$,
which satisfy the compatibility conditions 
\begin{equation}
  \label{eq:1.8}
  r_e f_{i(e)} = r_ef_{t(e)}
\end{equation}
for all edges, $e$ of $\Gamma$.

Following \cite{KR} we will define the $K$-theory analog of
$H_G(\Gamma)$ simply by replacing $\SS(\fg^*)$ in the definition by 
the \emph{representation ring}, $R(G)$, of $G$.

\end{definition}

\begin{definition}
$K_G(\Gamma)$ is the set of all functions, $f:V_{\Gamma} \to R(G)$, which 
satisfy the compatibility condition (\ref{eq:1.8}), $r_e$ being
the restriction map, $R(G) \to R(G_e)$.
\end{definition}

\begin{remarks*}
  \begin{enumerate}
  \item 
    Since $G$ is an $n$-torus, the representation ring $R(G)$ can be 
    identified with the \emph{character ring} of $G$, \emph{i.e.} the ring 
    of all finite sums
    \begin{equation}
      \label{eq:1.9}
      \sum m_k e^{2\pi i\alpha_k}
    \end{equation}
with $m_k \in \ZZ$ and $\alpha_k \in \ZZ^*_G$. We will frequently use this 
identification, referring to a representation by indicating the element of the 
character ring it corresponds to and vice-versa.

\item 
  Point-wise multiplication makes $K_G(\Gamma)$ into a
  ring.  Moreover, since the \emph{constant} functions satisfy
  (\ref{eq:1.8}), this ring contains the ring, $R(G)$, as a subring.

  \end{enumerate}
\end{remarks*}

Given $f \in K_G (\Gamma)$ let
\begin{equation}
  \label{eq:1.10}
  \chi(f) = \sum_{p \in V} f_p \prod_{i(e)=p}
  (1-e^{2\pi i \alpha_e})^{-1} \, .
\end{equation}
We will call $\chi (f)$ the \emph{Atiyah-Bott character} of the
class, $f$.  The individual summands on the right hand side are
elements of a quotient ring of $R(G)$; however, we will prove

\begin{theorem}
  \label{th:1.1}
  The sum (\ref{eq:1.10}) is an element of $R(G)$.  
\end{theorem}

Thus (\ref{eq:1.10}) defines a morphism of $R(G)$-modules
  \begin{displaymath}
    \chi : K_G (\Gamma) \to R(G)
  \end{displaymath}
which we will call the \emph{character map}.  A helpful way of
looking at this map is in terms of virtual representations.
Namely, to each $p \in V$, one can attach an infinite-dimensional 
virtual representation, $Q(\tau_p)$, the ``$\spin^{\CC}$-quantization'' 
of the action, $\tau_p$, of $G$ on $\CC^d$, and
Theorem \ref{th:1.1} asserts that the sum
\begin{equation}
  \label{eq:1.11}
  Q(f) = \bigoplus_{p \in V} Q (\tau_p) \otimes f_p
\end{equation}
is a finite dimensional virtual representation and that its
  character is given by (\ref{eq:1.10}).

Suppose, in particular, that $f$ has the form
\begin{equation}
  \label{eq:1.12}
  f_p = e^{2\pi i \alpha_p}, \quad \alpha_p \in \ZZ^*_G \, .
\end{equation}
Then by (\ref{eq:1.8})
\begin{equation}
  \label{eq:1.13}
  \alpha_q-\alpha_p=m_e\alpha_e
\end{equation}
for every pair of vertices, $p$ and $q$, and edge, $e$, joining $
p$ to $q$.

\begin{definition}
$f$ is symplectic if $m_e>0$ for all $e$.
\end{definition}

If $f$ is symplectic, the representation (\ref{eq:1.11}) has the
following convexity property.  (Compare with \cite[Theorem~6.3]{GS}.)

\begin{theorem}
  \label{th:1.2}
  If $\alpha$ is a weight of $Q(f)$ then $\alpha$ is in the
  convex hull of $\{ \alpha_p; \; p \in V \}$.
\end{theorem}

Let's denote this convex hull by $\Delta$.  We will call a
weight, $\alpha$, an \emph{extremal} weight if it is a vertex of
$\Delta$.  For these weights we will prove

\begin{theorem}
  \label{th:1.3}
  If $\alpha$ is extremal, it occurs in $Q(f)$ with
  multiplicity~$1$.
\end{theorem}

For non-extremal weights we will prove a more refined result. 
Fix a vector, $\xi$, in $\fg$ with the property $\alpha_e (\xi) \neq 0$ 
for all edges, $e$, of $\Gamma$; given a vertex, $p$, let 
$$\E_p = \{ e \in E_{\Gamma} ;  \; p \mbox{ is a vertex of } e  \mbox{ and }
\alpha_e(\xi) > 0 \},$$
and let $\sigma_p$ be the number of edges $e \in \E_p$ for which $p=t(e)$.
   
For $e \in \E_p$ define
$$ (-1)^e = \begin{cases} 
1 & \text{ if  } \;  p=i(e) \\ -1 & \text{ if  } \; p=t(e) \end{cases} \; ,$$ 
and let
\begin{eqnarray*}
(-1)^p & = & \prod_{e \in \E_p} (-1)^e = (-1)^{\sigma_p} \\
\delta_p & = & \tfrac{1}{2} \sum \alpha_e , \quad e \in \E_p \\
\delta^{\#}_p & = & \tfrac{1}{2} \sum (-1)^e \alpha_e
, \quad e \in \E_p \, .
\end{eqnarray*}

\begin{definition}
The Kostant partition function
\begin{displaymath}
  N_p : \ZZ^*_G \to \NN
\end{displaymath}
is the function which assigns to every weight, $\alpha$, the
number of distinct ways in which $\alpha$ can be written as a sum
\begin{displaymath}
  \alpha = \sum n_e\alpha_e , \quad e \in \E_p
\end{displaymath}
with non-negative integer coefficients.  
\end{definition}

\begin{theorem}
  \label{th:1.4}
The multiplicity with which a weight, $\alpha$, occurs in $Q(f)$
is equal to
\begin{equation}
  \label{eq:1.14}
  \sum_p (-1)^p N_p (\alpha -\alpha_p + \delta^{\#}_p
  - \delta_p) \, .
\end{equation}
(Compare with \cite[(1.13)]{GLS}.)
\end{theorem}

The next results which we will describe involve a
graph-theoretical analog of the notion of ``reduction by a circle 
action'' in symplectic geometry.  Let $T$ be a circle subgroup of 
 $G$ which is not contained in any of the groups, $G_e$.  Then if 
 $\xi$ is the infinitesimal generator of $T$
 \begin{equation}
   \label{eq:1.15}
   \alpha_e (\xi) \neq 0
 \end{equation}
for all $e$. A function, $\phi : V \to \RR$ is called \emph{a
  $T$-moment map} if for all edges $e \in E_{\Gamma}$
\begin{equation}
  \label{eq:1.16}
  \frac{\phi (t(e))-\phi(i(e))}
  {\alpha_e (\xi)} >0 \, .
\end{equation}
We recall (\cite[\S~2.2]{GZ1}) that there is a simple necessary
and sufficient condition for the existence of such a map.  By
(\ref{eq:1.6}) one can orient $\Gamma$ by assigning to each 
unoriented edge
the orientation for which $\alpha_e (\xi) >0$.  Then, for the
existence of a $T$-moment map, it is necessary and sufficient that 
this graph \emph{have no oriented cycles}.  We will call the
numbers, $\phi (p)$, the \emph{critical values} of $\phi$.  By
perturbing $\phi$ slightly one can arrange that these $\phi (p)$'s 
are all distinct.

Let $c \in \RR$ be a regular (non-critical) value of $\phi$, and
let $V_c$ be the set of all oriented edges, $e$, of $\Gamma$ with 
$\phi (t(e))>c>\phi (i(e))$.  One can make $V_c$ into the set of
vertices of a new object, $\Gamma_c$, and this object is our
graph-theoretical ``reduction of $\Gamma$ at $c$''.
(Unfortunately, $\Gamma_c$ is not a graph.  It is a slightly
more complicated object:  a ``hypergraph''.  For details see
\cite[\S~3]{GZ3}.)

Now fix an element, $f$, of $K_G (\Gamma)$, and for every edge, $
e$, in $V_c$ let $p=i(e)$ and let
\begin{equation}
  \label{eq:1.17}
  \hat{f}_e = f_p \prod_{e'} (1-e^{2\pi i \alpha_{e'}})^{-1}
\end{equation}
the product being over all $e'$ with $i(e')=p$ and $e'\neq e$.
By composing the inclusion map of $G_c$ into $G$ with the
projection of $G$ onto $G/T$, one gets a surjective finite-to-one 
map $\pi_e : G_e \to G/T$ and hence a ``push-forward'' in
$K$-theory (see Section~\ref{sec:rg})
\begin{displaymath}
  (\pi_e)_* : R (G_e) \to R (G/T) \, .
\end{displaymath}
This can be formally extended to elements of the quotient ring 
of $R(G_e)$ of the form \eqref{eq:1.17}, and by applying it to 
(\ref{eq:1.17}) one gets, for every vertex of $\Gamma_c$, an
element
\begin{equation}
  \label{eq:1.18}
  f^{\#}_c (e) = (\pi_e)_* r_e \hat{f}_e
\end{equation}
of a quotient ring of $R (G/T)$.

\begin{theorem}
  \label{th:1.5}
The sum
  \begin{displaymath}
    \chi_c (f) = \sum_{e \in V_c} f^{\#}_c (e)
  \end{displaymath}
is in $R(G/T)$.
\end{theorem}

We will prove this by proving a stronger result. Let
\begin{displaymath}
  f^{\#}_p = f_p \prod_{i(e)=p} (1-e^{2\pi i \alpha_e})^{-1}
\end{displaymath}
be the $p$\st{th} summand on the right hand side of
(\ref{eq:1.10}); and, for $g \in G$, consider the integral over $
T$
\begin{equation}
  \label{eq:1.19}
  \int f^{\#}_p (gt) \, dt \, .
\end{equation}
We will see in Section~\ref{sec:residue} that the integrand has poles at a
finite number of points, $t_i \in T$ so this integral as it
stands isn't well defined.  However, one can ``regularize'' it by 
moving the contour of integration to a curve in $T^{\CC}$ which
surrounds the $t_i$'s; and, denoting this regularized integral by
$\Res_T f^{\#}_p$ we will prove:

\begin{theorem}
  \label{th:1.6}
$\Res_T f^{\#}_p$ is an element of $R(G/T)$.

\end{theorem}

Our strengthened version of Theorem~\ref{th:1.5} asserts:

\begin{theorem}
  \label{th:1.7}
$\chi_c(f)$ is equal to the sum
\begin{equation}
  \label{eq:1.20}
  \sum_{\phi (p) >c} Res_T f_p^{\#}\, .
\end{equation}

\end{theorem}

Next we will explain what ``quantization commutes with
reduction'' translates into the context of graphs.  Recall that
an element, $f$, of $K_G(\Gamma)$ of the form
\begin{displaymath}
  f_p = e^{2\pi i \alpha_p}, \quad
  \alpha_p \in \ZZ^*_G
\end{displaymath}
is \emph{symplectic} if
\begin{displaymath}
  \alpha_q-\alpha_p = m_e \alpha_e , \quad m_e >0
\end{displaymath}
for every pair of vertices, $p$ and $q$, and edge, $e$ joining $
p$ to $q$.  For $f$ symplectic, the map
\begin{displaymath}
  \phi : V \to \RR, \quad p \to \alpha_p (\xi)
\end{displaymath}
is a $T$-moment map.  Assume zero is a regular value of this map, 
\emph{i.e.} $\alpha_p (\xi) \neq 0$ for all $p$; and let
$\Gamma_{red}=\Gamma_0$ and $\chi_{red}=\chi_0$ .

\begin{theorem}
  \label{th:1.8}
Let $Q (\Gamma)$ be the virtual representation of $G$ with
character, $\chi(f)$ and $Q(\Gamma_{red})$ the virtual
representation of $G/T$ with character, $\chi_{red}(f)$.  Then, as
virtual representations of $G/T$
\begin{equation}
  \label{eq:1.21}
  Q(\Gamma_{red}) = Q(\Gamma)^T \, .
\end{equation}

\end{theorem}

Finally in the last section of this paper we will show that if
$M$ is a GKM manifold and $\Gamma$ is its ``one-skeleton'', these
theorems about graphs have $K$-theoretic implications for $M$
(thanks to a beautiful recent result of Allen~Knutson and
Ioanid~Rosu which asserts that 
$K_G (M)\otimes \CC \simeq K_G(\Gamma) \otimes \CC$).

\section{Some algebraic preliminaries}
\label{sec:algprel}

We will collect in this section some elementary facts about lattices and 
tori which will be needed in the proofs. Let $V$ be an 
$n$-dimensional real vector space and let $L$ be a rank $n$ lattice 
sitting inside $V$. Let
$$L^* = \{ \alpha \in V^* ; \alpha(v) \in \ZZ \mbox{ for all } v \in L \}$$
be the dual lattice in $V^*$. An element of $L$ is \emph{primitive} if 
it is not of the form, $kv$, with $v \in L$ and $|k|>1$.
 
\begin{lemma}
\label{lem:2.1}
$v \in L$ is primitive if and only if there is an $\alpha \in L^*$ 
with $\alpha(v)=1$.
\end{lemma}

\begin{lemma}
\label{lem:2.2}
$v$ is primitive if and only if there exists a basis $v_1,...,v_n$ of 
$L$ with $v=v_1$.
\end{lemma}

Now let $G$ be an $n$-dimensional torus and let $\fg$ be its Lie algebra.

\begin{definition}
\label{def:2.3}
The \emph{group lattice} of $G$, $\ZZ_G$, is the kernel of the exponential 
map, 
$\exp~:~\fg~\to~G$ and its dual, $\ZZ_G^*$, is the \emph{weight lattice} of 
$G$.
\end{definition}

In particular
$$G=\fg / \ZZ_G$$
and the exponential map is just the projection of $\fg$ onto $\fg / \ZZ_G$.
Given a weight, $\alpha \in \ZZ_G^*$, let $\chi_{\alpha}$ be the character of
$G$ defined by
$$\chi_{\alpha}(g) = e^{2\pi i \alpha(x)}, \quad  g=\exp{(x)}.$$

\begin{proposition}
\label{prop:2.4}
If $\alpha$ is primitive, the subgroup
\begin{equation}
\label{eq:2.1}
G_{\alpha} = \{ g \in G ; \chi_{\alpha}(g)=1 \}
\end{equation}
is connected, i.e. is an $(n-1)$-dimensional subtorus of $G$. 
More generally, if $\beta$ is primitive and $\alpha=k\beta$, $k >1$, the 
identity component of $G_{\alpha}$ is $G_{\beta}$ and 
$G_{\alpha} / G_{\beta}$ is a finite cyclic group of order $k$.
\end{proposition}

Let $\xi$ be a primitive element of $\ZZ_G$ and let 
\begin{equation}
\label{eq:2.2}
G_{\xi} = \{ \exp{(t\xi)} ; 0 \leq t <1 \}.
\end{equation}
then $G_{\xi}$ is a closed, connected one-dimensional subgroup of $G$.

\begin{proposition}
\label{prop:2.5}
If $\alpha(\xi) =0$, then $G_{\xi} \subset G_{\alpha}$ and if 
$\alpha(\xi) \neq 0$, then $G_{\xi} \cap G_{\alpha}$ is a finite 
cyclic subgroup of order $|\alpha(\xi)|$.
\end{proposition}

Let $G_1=G/G_{\xi}$ and let $\gamma : G_{\alpha} \to G_1$ be the 
composition of the inclusion, $G_{\alpha} \to G$, and the projection, 
$G \to G_1$.

\begin{corollary}
\label{cor:2.5}
The map $\gamma$ is surjective and its kernel is a cyclic subgroup of 
$G_{\alpha}$ of order $|\alpha(\xi)|$.
\end{corollary}

\section{The representation ring $R(G)$}
\label{sec:rg}

The groups, $G$, in this section will be compact commutative Lie groups. 
For such a group every irreducible representation is one-dimensional, 
\emph{i.e.} is defined by a homomorphism of $G$ into $S^1$. Thus the 
elements of $R(G)$ can be identified with the \emph{character ring} 
of $G$: all finite sums of the form 
\begin{equation}
\label{eq:3.1}
\sum m_i \chi_i
\end{equation}
$m_i$ being an integer and $\chi_i$ a homomorphism of $G$ into $S^1$ 
(or ``character''.) Hence, if $G$ is an $n$-torus, 
\eqref{eq:3.1} is a sum of the form \eqref{eq:1.9}.

In this section we will discuss some functorial properties of this ring. 
First we note that $R(G)$ is naturally a contravariant functor, \emph{i.e.} 
if $\gamma : G \to H$ is a homomorphism of Lie groups, then a 
representation of $H$ can be converted, by composition with $\gamma$, 
into a representation of $G$; so there is a natural map
\begin{equation}
\label{eq:3.2}
\gamma^* : R(H) \to R(G)
\end{equation}
and it is easy to see that this is an algebra homomorphism. 
A much more interesting object for us will be a map in the opposite direction
\begin{equation}
\label{eq:3.3}
\gamma_* :R(G) \to R(H)
\end{equation}
which we will define here modulo the assumption
$$ (*) \qquad \mbox{\emph{ the kernel and cokernel of }} \gamma 
\mbox{\emph{ are finite.}}$$
First let's assume that $\gamma$ is surjective, \emph{i.e.} that $H=G/W$ 
and that $W$ is a finite subgroup of $G$. Let $\rho$ be a representation 
of $G$ on a vector space, $V$, and let $V^W$ be the vectors in $V$ which 
transform trivially under $W$. Then the restriction of $\rho$ to $V^W$ is 
a representation, $\rho^W$, of $G/W$ and $\gamma_*$ is the map defined by 
$\rho \to \rho^W$.

Next assume that $\gamma$ is injective, \emph{i.e.} that $G$ is a closed 
subgroup of $H$ and $G \backslash H$ is finite. Given a representation, 
$\rho$, of $G$ on a vector space, $V$, let $\rho_{ind}$ be the induced 
representation of $H$ (\emph{i.e.} let $V_{ind}$ be the vector space 
consisting of maps $f : H \to V$ which satisfy $f(gh) = \rho(g)f(h)$ and 
let 
$$(\rho_{ind} f)(k) = f(kh^{-1})$$
for all $k \in H$). In this case $\gamma_*$ is the map defined by 
$\rho \to \rho_{ind}$.

Finally if $\gamma$ is neither injective nor surjective, let $G_1$ be the 
image of $\gamma$. Then $\gamma$ factors into the submersion 
$\gamma_1 : G \to G_1$, composed with the inclusion, $\gamma_2 : G_1 \to H$, 
and one defines
\begin{equation}
\label{eq:3.5}
\gamma_* = (\gamma_2)_* (\gamma_1)_*.
\end{equation}
This map is unfortunately not a ring homomorphism, but it is a 
morphism of $R(H)$-modules: for $\chi \in R(G)$ and $\tau \in R(H)$
\begin{equation}
\label{eq:3.6}
\gamma_*(\chi \gamma^* \tau) = (\gamma_* \chi) \tau.
\end{equation}

We will mostly be interested in the case when $\gamma$ is a submersion, 
\emph{i.e.} when $H = G/W$. In this case one has an alternative way of 
looking at $\gamma$:

\begin{lemma}
\label{lem:3.1}
Let $\rho$ be a unitary representation of $G$ on a complex vector space, $V$. 
Then the orthogonal projection of $V$ onto $V^W$ is given by the operator
\begin{equation}
\label{eq:3.7}
P = \frac{1}{|W|} \sum_{w \in W} \rho(w).
\end{equation}
\end{lemma}

\begin{proof}
If $v \in V^W$, $\rho(w)v=v$, so $Pv=v$. Moreover, for all 
$v \in V$ and $a \in W$,
$$\rho(a) Pv = \frac{1}{|W|} \sum_{w\in W} \rho(aw) v = Pv,$$
so $Pv \in V^W$. Finally, 
$$P^* = \frac{1}{|W|} \sum_{w\in W} \rho(w^{-1}) =P,$$
so $P$ is the \emph{orthogonal} projection of $V$ onto $V^W$.
\end{proof}

\begin{corollary}
\label{cor:3.2}
Let $g$ be an element of $G$ and let $\bar{g}$ be its image in $G/W$. Then
\begin{equation}
\label{eq:3.8}
(\gamma_* \rho)(\bar{g}) = \frac{1}{|W|} \sum_{w \in W} \rho (gw).
\end{equation}
In particular, let $f : G \to \CC$ be the function \eqref{eq:3.1}, 
\emph{i.e.} 
\begin{equation}
\label{eq:3.9}
f(g)= \sum m_i \chi_i(g).
\end{equation}
Then
\begin{equation}
\label{eq:3.10}
\gamma_*f(\bar{g}) = \frac{1}{|W|} \sum_{w \in W} f(gw).
\end{equation}
\end{corollary}

\section{Convexity and multiplicities}
\label{sec:proof1.1}

\begin{proof}[The proof of Theorem~\ref{th:1.1}]: \\
Let $\a_1, \cdots , \a_N$ be primitive vectors such that 
for every $e \in E_{\Gamma}$ there exists a unique
$k \in \{1,...,N\}$ such that $\a_e$ is a multiple of $\a_k$.
If $m_1\a_1,...,m_s\a_1$ are all the occurrences of multiples of 
$\a_1$ among all the weights, let 
$M_1 = \mbox{l.c.m.}( m_1,...,m_s ). $
Similarly we define $M_2,...,M_N$. Then 
\begin{equation}
   \label{eq:2}
   \chi ( f) = \frac{g}{\prod_{j=1}^N (1-\e{M_j\alpha_j})}
\end{equation}
with $g \in R(G)$. We will show that $1-\e{M_1\a_1}$ divides $g$
in $R(G)$.

The vertices of $\Gamma$ can be divided into two categories:

\begin{enumerate}
\item The first subset, $V_1$, contains the vertices 
$p \in V_{\Gamma}$ for which none of the $\alpha_{e}$'s with $i(e)=p$, 
is a multiple of $\alpha_1$

\item The second subset, $V_2$, contains the vertices 
$p \in V_{\Gamma}$ for which there exists an edge $e$ such that
$i(e)=p$ and $\alpha_{e}$ is a multiple of $\alpha_1$. (Notice that 
there will be exactly one such edge.)

\end{enumerate}

The part of \eqref{eq:1.10} corresponding to vertices in the first 
category will then be of the form
\begin{equation}
	\label{eq:3}
    \sum_{p \in V_1} f_p \prod_{i(e)=p} (1-\e{\alpha_{e}})^{-1} =
    g_1 \prod_{j=2}^N (1-\e{M_j \alpha_j})^{-1} 
\end{equation}
with $g_1 \in R(G).$

If $p \in V_2$ then there exists an edge $e$ issuing from $p$ such that 
$\alpha_e = m \alpha_1$ with 
$m \in \ZZ - \{ 0 \}$; let $q= t(e)$.
Since $\alpha_{\bar{e}} = - \alpha_{e}$
it follows that $q \in V_2$ as well and thus the vertices in $V_2$
can be paired as above.

Let $e_k, k=1, \ldots,d$ and $e'_k , k=1, \ldots ,d$ be the edges
issuing from $p$ and $q$ respectively, with $e_d=e$, $e'=\bar{e}$.  
Then, by \eqref{eq:1.4}, the $e_k$'s can be ordered so that
$$ r_e (\e{\a_{e_k}}) = r_e (\e{\a_{e'_k}}) \; , 
\mbox{ for all } k=1,...,d-1 $$
which implies that
\begin{equation}
  \label{eq:4}
   1- \e{\a_{e_k}} \equiv  1- \e{\a_{e'_k}} 
            \pmod{ 1- \e{\alpha_e}} .
\end{equation}
Similarly, from
$$r_e (f_p) = r_e ( f_q) $$
we deduce that
\begin{equation}
   \label{eq:5}
  f_q \equiv f_p \pmod{ 1-\e{\alpha_e}} \; .
\end{equation}
The part of \eqref{eq:1.10} corresponding to $p$ and $q$,
$$
 f_p\prod_{j=1}^d (1-\e{\alpha_{e_j}})^{-1}  + 
 f_q\prod_{j=1}^d (1-\e{\alpha_{e'_j}})^{-1} \; ,$$
can be expressed as
\begin{equation}
\label{eq:6}
\frac{f_p \prod_{j=2}^N (1-\e{\a_{e'_j}}) -
\e{m\a_1} f_q \prod_{j=2}^N (1-\e{\a_{e_j}}) }
{(1-\e{m\a_1}) \prod_{j=2}^N (1-\e{\a_{e_j}}) 
\prod_{j=2}^N (1-\e{\a_{e'_j}})}.
\end{equation}
From the congruences \eqref{eq:4} and 
\eqref{eq:5} we conclude that $1-\e{m\alpha_1}$ divides the numerator of
\eqref{eq:6}, so we deduce that
$$ 
f_p \prod_{j=1}^d (1-\e{\alpha_{e_j}})^{-1}  + 
 f_q\prod_{j=1}^d (1-\e{\alpha_{e'_j}})^{-1} = 
\frac{g_{p,q}}{\prod_{j=2}^N (1-\e{M_j\a_j})}
$$
with $g_{p,q} \in R(G)$. Therefore
\begin{equation}
  \label{eq:8}
    \sum_{p \in V_2} f(p)\prod_{i(e)=p} (1-\e{\alpha_{e}})^{-1} =
    \frac{g_2}{\prod_{j=2}^N (1-\e{M_j\alpha_j}) }
\end{equation}
with $g_2 \in R(G)$. Adding \eqref{eq:3} and \eqref{eq:8}
we obtain
$$
\frac{g}{\prod_{j=1}^N (1-\e{M_j\alpha_j})} =
\frac{ g_1 + g_2}{\prod_{j=2}^N (1-\e{M_j\alpha_j})}
$$
with $g_1 + g_2 \in R(G)$, hence $1-\e{M_1\alpha_1}$ divides $g$.
The same argument can be used to show that each $1-\e{M_j\alpha_j}$ 
divides $g$.

The proof of the theorem now follows from:

\begin{lemma}
 \label{lem:1}
If $P \in R(G)$ and $\a$, $\b$ are linearly independent weights such that
$1-\e{\a}$ divides $(1-\e{\b})P$, then $1-\e{\a}$ divides  $P$.
\end{lemma}

\begin{lemma}
\label{lem:2}
If $P \in R(G)$ and $\b_1,...,\b_k$ are pairwise linearly independent 
weights such that $1-\e{\b_j} \mbox{ divides } P$ for all $j=1,...,k$ then
$$ \qquad \qquad (1-\e{\b_1}) \cdots (1-\e{\b_k}) \mbox{ divides } P \; . 
\qquad \qquad \qquad  \qed$$
\end{lemma}
\renewcommand{\qed}{}
\end{proof}

\begin{proof}[The proof of Theorem~\ref{th:1.2}]: \\
Let $\alpha$ be a weight that is not in the convex hull of 
$\{ \alpha_p; p \in V_{\Gamma} \}$. Then there exists $\xi \in \fg$ and 
$p_0 \in V_{\Gamma}$ such that  $ (\alpha -\alpha_{p_0})(\xi) < 0 $
and $ (\alpha_p - \alpha_{p_0})(\xi) > 0$ for all $p \neq p_0$.
If $e \in E_{\Gamma}$ and $\alpha_e(\xi) <0$ then
$$(1-e^{2\pi i \alpha_e})^{-1} = - e^{2 \pi i \alpha_e} 
(1-e^{2\pi i \alpha_{\bar{e}}})^{-1},$$
and using this we deduce that
\begin{equation}
\label{eq:9.1}
\chi(f) = \sum_{p \in V} (-1)^{\sigma_p} 
e^{2 \pi i (\alpha_p - \sum' \alpha_e)}
\prod_{e \in \E_p} (1-e^{2 \pi i \alpha_e})^{-1},
\end{equation}
where $\sum' \alpha_e$ in the exponent is the sum
\begin{equation}
\label{eq:9.2}
\sum_{\substack{i(e)=p \\ \alpha_e(\xi) <0}} \alpha_e =
- \sum_{\substack{t(e)=p \\ \alpha_e(\xi) > 0}} \alpha_e = 
\delta_p^{\#} - \delta_p.
\end{equation}
From \eqref{eq:9.1} and \eqref{eq:9.2} we deduce that
\begin{equation}
\label{eq:9.3}
\chi(f) = \sum_{p \in V} (-1)^p 
e^{2 \pi i (\alpha_p - \sum' \alpha_e)}
\prod_{e \in \E_p} ( \sum_{k_e \geq 0} e^{2 \pi i k_e \alpha_e}).
\end{equation}

Suppose $\alpha$ is a weight of $Q(f)$; 
then there exists $p \in V_{\Gamma}$ and 
non-negative integers $\{k_e\}_{e \in \E_p}$ such that
\begin{equation}
\label{eq:9.4}
\alpha = \alpha_p + \sum_{\substack{t(e)=p \\ \alpha_e(\xi) > 0}} \alpha_e +
\sum_{e \in \E_p} k_e \alpha_e,
\end{equation}
which implies that
\begin{equation}
\label{eq:9.5}
\alpha - \alpha_{p_0}= \alpha_p -\alpha_{p_0} + 
\sum_{\substack{t(e)=p \\ \alpha_e(\xi) > 0}} \alpha_e +
\sum_{e \in \E_p} k_e \alpha_e.
\end{equation}

But when we evaluate \eqref{eq:9.5} at $\xi$, the right hand side is 
non-negative, while the left hand side is strictly negative ! This 
contradiction proves that $\alpha$ is not a weight of $Q(f)$.
\end{proof}

\begin{proof}[The proof of Theorem~\ref{th:1.3}]: \\
Let $\alpha=\alpha_{p_0}$ be an extremal weight, 
\emph{i.e.} a vertex of $\Delta$. Then there exists $\xi \in \fg$ such that
$(\alpha_p - \alpha_{p_0}) (\xi) > 0$ for all $p \neq p_0$. 
In this case \eqref{eq:9.5} implies
\begin{equation}
\label{eq:9.6}
0=  (\alpha_p -\alpha_{p_0})(\xi) + 
\sum_{\substack{t(e)=p \\ \alpha_e(\xi) > 0}} \alpha_e(\xi) +
\sum_{e \in \E_p} k_e \alpha_e(\xi) . 
\end{equation}
Since each term on the right hand side is non-negative, \eqref{eq:9.6} 
is only true if  
\begin{enumerate}
\item $p=p_0$ (which also implies that $\alpha_e(\xi) <0$ for all 
$e$ with $t(e)=p$, \emph{i.e.} that there are no terms in the first sum), and
\item $k_e=0$ for all $e \in \E_p$.
\end{enumerate}

This proves that the multiplicity with which $\alpha$ occurs in $Q(f)$ is 1.
\end{proof}

\begin{proof}[The proof of Theorem~\ref{th:1.4}]: \\
From \eqref{eq:9.2} and 
\eqref{eq:9.4}, the multiplicity with which a weight $\alpha$ appears in 
the term corresponding to the vertex $p$ is equal to $(-1)^p$ times the 
number of distinct ways in which $\alpha - \alpha_p + \delta_p^{\#}-\delta_p$ 
can be written as a sum
$$\sum_{e \in \E_p} k_e\alpha_e,$$
with $k_e$'s non-negative integers; and this number is 
$N_p(\alpha - \alpha_p + \delta_p^{\#}-\delta_p)$.
Counting the contributions given by all the vertices we 
obtain \eqref{eq:1.14}.
\end{proof}

\section{The residue operation}
\label{sec:residue}

Let $G$ be an $n$-dimensional torus, let $T$ be a circle subgroup of $G$, 
and let 
$$\chi_k= e^{2 \pi i \alpha_k}, \quad k=1,...,d$$
be characters of $G$ and $f$ an element of the character ring, $R(G)$. 
The goal of this section is to make sense of the integral
\begin{equation}
\label{eq:4.1}
\int_T \frac{f(gt)}{\prod(1-\chi_k(gt))} \; dt
\end{equation}
as a function of $g \in G$. If the restriction of $\chi_k$ to $T$ is 
identically one, the denominator in the integrand is identically zero 
when $\chi_k(g)=1$. Hence, for \eqref{eq:4.1} to make sense, we are 
forced to assume that the restriction of $\chi_k$ to $T$ is \emph{not} 
identically one. Even with this assumption, however, the integrand has poles 
at the points where $\chi_{k}(gt)=1$; so to make sense of \eqref{eq:4.1} we 
must ``regularize'' this integral and this we will do as follows. Fix a 
basis vector, $\xi$ of $\ZZ_T$, and identify $T$ with $S^1$ via the map 
$$\exp{(s\xi)} \to e^{2\pi i s}.$$
Then, with $z=e^{2\pi i x}$, the integrand of \eqref{eq:4.1} becomes a 
meromorphic function
\begin{equation}
\label{eq:4.2}
f^{\#}(gz) = f(gz)\prod_{k=1}^d (1-\chi_k(gz))^{-1}
\end{equation}
on the complex plane with poles on the unit circle. Now move the contour 
of integration from the unit circle to a contour surrounding these poles, 
\emph{e.g.} a circle of radius greater than one oriented in a 
\emph{counter-clock-wise} sense plus a circle of radius less than 
one oriented in a \emph{clock-wise} sense. In other words, 
replace \eqref{eq:4.1} by the integral
\begin{equation}
\label{eq:4.3}
\frac{1}{2\pi i}\int_{C_+} f^{\#} (gz) \frac{dz}{z}
\end{equation}
minus the integral
\begin{equation}
\label{eq:4.4}
\frac{1}{2\pi i}\int_{C_-} f^{\#} (gz) \frac{dz}{z} ,
\end{equation}
$C_+$ being a circle of radius greater than one and $C_-$ a circle of 
radius less than one, both these circle being oriented in a 
counter-clock-wise sense. Let us denote this regularized integral, 
\emph{i.e.} the difference of \eqref{eq:4.3} and \eqref{eq:4.4}, by 
$(Res_T f^{\#})(g)$. It is easy to see that this function is $T$-invariant,
$$(Res_T f^{\#})(gt) = (Res_T f^{\#})(g)$$
and hence defines a function on $G/T$. We will prove:

\begin{theorem}
\label{th:4.1}
$Res_T f^{\#}$ is an element of $R(G/T)$.
\end{theorem}

\begin{remark*}
The definition of $Res_T f^{\#}$ depends on the identification of $S^1$ 
with $T$ given by $\exp{(x\xi)} \leftrightarrow e^{2\pi i x}$. If we 
replace $\xi$ by $-\xi$, the orientations of the circles, $C_+$ and $C_-$, 
will get reversed, and hence this will change the signs of \eqref{eq:4.3} 
and \eqref{eq:4.4}.
\end{remark*}

In proving Theorem \ref{th:4.1}, we can assume without loss of generality that
$f=e^{2 \pi i \alpha}$, $\alpha \in \ZZ_G^*$. Let $e_1,..,e_n$ be a basis of 
$\ZZ_G$ with $\xi=e_n$ and let $y_1,...,y_{n-1}$ and $x$ be the coordinates on 
$\fg$ associated with this basis. We can then write 
\begin{equation}
\label{eq:4.5}
\alpha_i = k_i x + \beta_i(y)
\end{equation}
and
\begin{equation}
\label{eq:4.6}
\alpha = kx + \beta(y)
\end{equation}
with $k_i = \alpha_i(\xi)$ and $k=\alpha(\xi)$. Thus letting
\begin{eqnarray}
\label{eq:4.7}
z & = & e^{2 \pi i x} \\
\label{eq:4.8}
a_i & = & e^{2 \pi i \beta_i(y)} \\
\label{eq:4.9}
b & = & e^{2 \pi i \beta(y)}
\end{eqnarray}
the integrals \eqref{eq:4.3} and \eqref{eq:4.4} become 
\begin{equation}
\label{eq:4.10}
\frac{1}{2\pi i} \int_{C_+} \frac{bz^k}{\prod (1-a_iz^{k_i})} \frac{dz}{z}
\end{equation}
and
\begin{equation}
\label{eq:4.11}
\frac{1}{2\pi i} \int_{C_-} \frac{bz^k}{\prod (1-a_iz^{k_i})} \frac{dz}{z}.
\end{equation}
Therefore, to prove the theorem, it suffices to show that each of these 
integrals individually is in $R(G/T)$. To verify this for \eqref{eq:4.11},
let's order the factors in the denominator of the integrand so that 
$k_i = - k_i' < 0$ for $1 \leq i \leq r$ and $k_i > 0$ for 
$r+1 \leq i \leq d$. This integrand is then equal to
\begin{equation}
\label{eq:4.12}
\frac{bz^{k'-1}}{\prod_{i=1}^r(z^{k_i'} - a_i) \prod_{i=r+1}^d (1-a_iz^{k_i})}
\end{equation}
with $k'=k-k_1 - ... - k_r$. Hence, if $k'>0$, \eqref{eq:4.12} is 
\emph{holomorphic} at zero; so, in particular:

\begin{lemma}
\label{lem:4.2}
The integral \eqref{eq:4.4} is zero if $k > k_1 + ... + k_r$.
\end{lemma}

For $1 \leq i \leq r$ and $z \approx 0$, let $a_i' = a_i^{-1}$ and let 
\begin{equation}
\label{eq:4.13}
S_i(z) = \frac{1}{z^{k_i'}-a_i} = \frac{-a_i}{1-a_i'z^{k_i'}} =
-a_i \sum_{l=0}^{\infty} (a_i'z^{k_i'})^l,
\end{equation}
and, for $r+1 \leq i \leq d$, let
\begin{equation}
\label{eq:4.14}
S_i(z) = \frac{1}{1-a_iz^{k_i}} = \sum_{l=0}^{\infty} (a_iz^{k_i})^l.
\end{equation}
Then the integral \eqref{eq:4.11} is just the degree -1 term in the 
Laurent series
\begin{equation}
\label{eq:4.15}
bz^{k'-1} \prod_{i=1}^d S_i(z)
\end{equation}
and this term is clearly a polynomial in 
$b, a_1,..,a_r, a_1^{-1},...,a_r^{-1}$, and $a_{r+1},...,a_d$, with 
integer coefficients. Hence, by \eqref{eq:4.8} and \eqref{eq:4.9}, it is 
clearly a trigonometric polynomial in $y_1,..,y_{n-1}$. From this, 
together with Lemma \ref{lem:4.2}, we conclude:

\begin{theorem}
\label{th:4.3}
The integral \eqref{eq:4.4} is an element of $R(G/T)$. Moreover, if 
$f=e^{2 \pi i \alpha}$ and $k = \alpha(\xi)$, this integral is zero if 
$k > k_1 + ... +k_r$.
\end{theorem}

To evaluate the integral \eqref{eq:4.10} we make the substitution, 
$z \to z^{-1}$ and reduce this integral to an integral of the type we've just
evaluated. We conclude

\begin{theorem}
\label{th:4.4}
The integral \eqref{eq:4.3} is an element of $R(G/T)$. Moreover, if 
$f=e^{2 \pi i \alpha}$ and $k = \alpha(\xi)$, this integral is zero if 
$k < k_{r+1} + ... + k_d$.
\end{theorem}

Since $k_i$ is negative for $1 \leq i \leq r$ and positive for 
$r+1 \leq i \leq d$ we have, in particular:

\begin{proposition}
\label{prop:4.5}
If $f=e^{2\pi i \alpha}$ with $k=\alpha(\xi)$, the integral \eqref{eq:4.4}
is zero when $k$ is positive and the integral \eqref{eq:4.3} is zero when 
$k$ is negative. Moreover, if $k=0$, \eqref{eq:4.4} is zero when $r >0$
and \eqref{eq:4.3} is zero when $d-r > 0$.
\end{proposition}

Suppose now that the weights, $\alpha_i$, $i=1,..,d$ are 
\emph{pairwise linearly independent}, \emph{i.e.} suppose that $\alpha_i$ 
and $\alpha_j$ are linearly independent for $i \neq j$. Then the integrand in 
\eqref{eq:4.10} - \eqref{eq:4.11}:
\begin{equation}
\label{eq:4.16}
bz^{k-1}\prod_{i=1}^d (1-a_iz^{k_i})^{-1}
\end{equation}
has simple poles on the unit circle for generic values of $y$. 
(Recall that since $a_k=e^{2\pi i \beta_k(y)}$, the location of these 
poles depends on $y$.) Thus, one can compute the difference between 
\eqref{eq:4.10} and \eqref{eq:4.11} by computing the residues of 
\eqref{eq:4.16} at these poles. We will 
show that the sum of these residues, which is, 
by definition, the regularized integral \eqref{eq:4.1}, is given by an 
expression involving the $K$-theoretic push-forward which we described 
in Section \ref{sec:rg}. More explicitly, let $G_i$ be the kernel of 
the homomorphism $\chi_i : G \to S^1$ and let $r_i$ be the restriction map 
$R(G) \to R(G_i)$ and $\pi_i$ the projection of $G_i$ onto $G/T$. We will 
prove:

\begin{theorem}
\label{th:4.6}
Let $f^{\#}$ be the function \eqref{eq:4.2} and let 
$$\hat{f_i} = f \prod_{j \neq i} (1 - \chi_j)^{-1}.$$
Then 
\begin{equation}
\label{eq:4.17}
Res_Tf^{\#} = \sum_{i=1}^r (\pi_i)_* r_i \hat{f_i} -
\sum_{i=r+1}^d (\pi_i)_* r_i \hat{f_i}.
\end{equation}
\end{theorem}

%\section{The proof of Theorem \ref{th:4.6}}
%\label{sec:proof}

\begin{proof}
Let $\theta_1, ..., \theta_d$ be real numbers, let $k_1,...,k_d$ be 
integers and let $a_i=e^{2 \pi i \theta_i}$. As above 
we will order the $k_i$'s so that $k_i <0$ for $1 \leq i \leq r$ and 
$k_i > 0$ for $r+1 \leq i \leq d$. Let $g(z)$ be the function 
\eqref{eq:4.16}, \emph{i.e.}
$$g(z)= bz^{k-1} \prod_{i=1}^d (1-a_iz^{k_i})^{-1}.$$

\begin{lemma}
\label{lem:5.1}
Suppose that, for $i \neq j$, $\theta_i, \theta_j$ and 1 are 
linearly independent over the rationals. Then $g(z)$ has simple 
poles on the unit circle.
\end{lemma}

\begin{proof}
Let 
$$\omega_i = e^{2 \pi i/k_i} \quad \mbox{ and } \quad 
a_i^{-1/k_i}= e^{-2\pi i \theta_i/k_i}.$$
Then these poles are at the points
\begin{equation}
\label{eq:5.1}
\omega_i^l a_i^{-1/k_i}, \quad 1 \leq l \leq k_i, \; 1 \leq i \leq d;
\end{equation}
so if $\theta_i, \theta_j$ and 1 are linearly independent over the 
rationals these poles are distinct.
\end{proof}

Let us compute the residue of $g(z)$ at the pole \eqref{eq:5.1}. The 
quotient
$$\frac{z-\omega_i^l a_i^{-1/k_i} }{1-a_iz^{k_i}}$$
evaluated at $z=\omega_i^l a_i^{-1/k_i}$ is equal, by l'Hopital's rule,
to: 
$$\frac{1}{-a_ik_iz^{k_i-1}} \quad \mbox{ or, alternatively }
\quad \frac{z}{-a_ik_iz^{k_i}}$$
evaluated at $z=\omega_i^l a_i^{-1/k_i}$, and since 
$(\omega_i^l a_i^{-1/k_i})^{k_i}= a_i^{-1}$, this quotient is just
\begin{equation}
\label{eq:5.2}
-\frac{1}{k_i} \omega_i^la_i^{-1/k_i}.
\end{equation}

Thus the residue at $z=\omega_i^l a_i^{-1/k_i}$ of the function
$$g(z) = \frac{1}{1-a_iz^{k_i}} bz^{k-1} \prod_{j \neq i} 
(1-a_jz^{k_j})^{-1}$$
is just
$$-\frac{b}{k_i} (\omega_i^l a_i^{-1/k_i})^k \prod_{j \neq i}
(1-a_j( a_i^{-1/k_i} )^{k_j})^{-1}$$
which, if we set
\begin{equation}
\label{eq:5.3}
b_i=ba_i^{-k/k_i}
\end{equation}
and
\begin{equation}
\label{eq:5.4}
a_{j,i} = a_ja_i^{-k_j/k_i},
\end{equation}
can be written
\begin{equation}
\label{eq:5.5}
-\frac{1}{k_i} (\omega_i^l)^k b_i \prod_{j \neq i} 
(1-(\omega_i^l)^{k_j} a_{j,i})^{-1}.
\end{equation}

We will now show that if we give $b$ and $a_i$ the values 
\eqref{eq:4.8} - \eqref{eq:4.9} the sum of these residues is identical 
with the right 
hand side of \eqref{eq:4.17}. If $b$ is equal to \eqref{eq:4.9} and $a_i$ is 
equal to \eqref{eq:4.8}, then by \eqref{eq:4.5} and \eqref{eq:4.6}
\begin{equation}
\label{eq:5.6}
b_i = e^{2 \pi i \sigma_i}
\end{equation}
and
\begin{equation}
\label{eq:5.7}
a_{j,i} = e^{2 \pi i \alpha_{j,i}},
\end{equation}
where 
\begin{equation}
\label{eq:5.8}
\sigma_i = \alpha -\frac{k}{k_i}\alpha_i
\end{equation}
and
\begin{equation}
\label{eq:5.9}
\alpha_{j,i} = \alpha_j - \frac{k_j}{k_i} \alpha_i.
\end{equation}

Let's now give a more ``intrinsic'' definition of $\sigma_i$ and 
$\alpha_{j,i}$: Let $\fg_i$ be the Lie algebra of the group, $G_i$, and 
$\ft$ the Lie algebra of $T$. Since $G_i$ is by definition the kernel of 
the homomorphism, $e^{2 \pi i \alpha_i} : G \to S^1$, $\fg_i$ is the 
annihilator of $\alpha_i$; so, by \eqref{eq:5.8}, $\sigma_i$ is 
the \emph{unique} element of $\fg^*$ which is annihilated by $\ft$ and 
has the same restriction to $\fg_i$ as $\alpha$. Similarly, $\alpha_{j,i}$
is the unique element of $\fg^*$ which is annihilated by $\ft$ and has the 
same restriction to $\fg_i$ as $\alpha_j$. Note, by the way, that since 
$\sigma_i$ and $\alpha_{j,i}$ are annihilated by $\ft$, they are in the 
dual vector space to $\fg/\ft$; or, in other words, in the dual of the 
Lie algebra of $G/T$.

Consider the kernel of the map $G_i \to G/T$. This consists of the elements
$$ \exp{(\frac{l}{k_i}\xi)} , \quad l=1,..,k$$
and by \eqref{eq:4.5} and \eqref{eq:4.6}
\begin{equation}
\label{eq:5.10}
e^{2 \pi i \alpha} (\exp{(\frac{l}{k_i}\xi)}) = (\omega_i^l)^k
\end{equation}
and
\begin{equation}
\label{eq:5.11}
e^{2 \pi i \alpha_j} (\exp{(\frac{l}{k_i}\xi)}) = (\omega_i^l)^{k_j}.
\end{equation}

Thus the sum
$$-\frac{1}{k_i} \sum_{l=1}^{k_i} (\omega_i^l)^k e^{2 \pi i \sigma_i}
\prod_{j \neq i} (1 -(\omega_i^l)^{k_j} e^{2 \pi i \alpha_{j,i}})^{-1}$$
of the residues of $g(z)$ over the poles $(\omega_i^l)a_i^{-1/{k_i}}$, 
$1 \leq l \leq k_i$ is by formula \eqref{eq:3.10} identical to the expression
$$-(\pi_i)_* r_i \frac{ e^{2 \pi i \alpha}}
{\prod_{j \neq i} (1- e^{2 \pi i \alpha_j})}$$
if $r+1 \leq i \leq d$ (in which case $k_i = |k_i|$) and is equal to 
$$(\pi_i)_* r_i \frac{ e^{2 \pi i \alpha}}
{\prod_{j \neq i} (1- e^{2 \pi i \alpha_j})}$$
when $1 \leq i \leq r$ (in which case $k_i = -|k_i|$).
\end{proof}

\section{Quantization commutes with reduction}
\label{sec:quantred}

We will prove below Theorems \ref{th:1.7} and \ref{th:1.8} of
Section \ref{sec:intro}. As in Theorem \ref{th:1.7}, let $f$ be an 
element of $K_G(\Gamma)$, let $\phi : V \to \RR$ be a $T$-moment map, 
let $c$ be a regular value of $\phi$ and let $e$ be an oriented 
edge of $\Gamma$ with 
$ \phi(q) > c > \phi(p), $
where $p=i(e)$ and $q=t(e)$ 
(\emph{i.e.} $e$ corresponds to a vertex of the hypergraph, $\Gamma_c$; 
we will denote this vertex by $e$, as well.)

Consider the expressions 
\begin{eqnarray}
\label{eq:6.1}
\hat{f_e} & = & f_p \prod_{e'} (1 - e^{2 \pi i \alpha_{e'}})^{-1} \\
\label{eq:6.2}
\hat{f_{\bar{e}}} & = & f_q \prod_{e''} (1 - e^{2 \pi i \alpha_{e''}})^{-1}
\end{eqnarray}
the product in \eqref{eq:6.1} being over all edges, $e' \neq e$, with 
$i(e')=p$, and the product in \eqref{eq:6.2} being over all edges, 
$e'' \neq \bar{e}$, with $i(e'')=q$.

\begin{lemma}
\label{lem:6.1}
Let $r_e=r_{\bar{e}}$ be the restriction map $R(G) \to R(G_e)$. Then 
\begin{equation}
\label{eq:6.3}
r_e \hat{f_e} = r_{\bar{e}} \hat{f_{\bar{e}}}.
\end{equation}
\end{lemma}

Let $\pi_e = \pi_{\bar{e}}$ be the projection of $G_e$ onto $G/T$. 
As a corollary of Lemma \ref{lem:6.1} we get two alternative ways of 
defining \eqref{eq:1.18}:
\begin{equation}
\label{eq:6.4}
(\pi_e)_* r_e \hat{f_e} = (\pi_{\bar{e}})_* r_{\bar{e}} \hat{f_{\bar{e}}} =
f_c^{\#}(e),
\end{equation}
and, as a consequence of \eqref{eq:6.4}, the following theorem:

\begin{theorem}
\label{th:6.2}
Let $c$ and $c'$ be regular values of $\phi$. Suppose there exists 
just one vertex, $p$, with $c < \phi(p) < c'$. Then
\begin{equation}
\label{eq:6.5}
\chi_c(f) -\chi_{c'}(f) = Res_T 
\Bigl( f_p \prod_{i(e)=p} (1-e^{2 \pi i \alpha_e})^{-1} \Bigr)
\end{equation}
\end{theorem}

\begin{proof}
If $e \in V_c$ and $t(e) \neq p$, then $e \in V_{c'}$, 
and if $e \in V_{c'}$ and $i(e) \neq p$, 
then $e \in V_c$. Moreover, in both cases,
\begin{equation}
\label{eq:6.6}
f_c^{\#}(e) = f_{c'}^{\#}(e),
\end{equation}
by \eqref{eq:6.4}. Thus, if $e_i$, $i=1,..,r$ are the elements of 
$V_c$ with $t(e_i)=p$, and $e_i$, $i=r+1,...,d$, are the elements of 
$V_{c'}$ with $i(e_i)=p$, 
the difference between $\chi_c(f)$ and $\chi_{c'}(f)$ is, by
\eqref{eq:6.4}, equal to
$$\sum_{i=1}^r f_c^{\#}(e_i) - \sum_{i=r+1}^d f_{c'}^{\#}(e_i),$$
or, also by \eqref{eq:6.4}, to
$$\sum_{i=1}^r (\pi_{e_i})_* r_{e_i} \hat{f}_{e_i} -
\sum_{i=r+1}^d (\pi_{e_i})_* r_{e_i} \hat{f}_{e_i},$$
which, by \eqref{eq:4.17}, is identical with 
$$Res_T \Bigl( f_p \prod_{i(e)=p} (1-e^{2 \pi i \alpha_e})^{-1} \Bigr). 
\qquad \qquad \qquad \qed$$
\renewcommand{\qed}{}
\end{proof}

To prove Theorem \ref{th:1.7}, let $c_0 < c_1 < ... < c_N$ be regular 
values of $\phi$ with $c_0 =c$, $c_N$ greater than $\phi_{max}$, and 
with only one critical point, $p_i$, between $c_i$ and $c_{i+1}$. Then
$$ \chi_c(f) = \sum _{i=0}^N (\chi_{c_i}(f) - \chi_{c_{i+1}}(f)) =
\sum_{\phi(p) > c} 
Res_T \Bigl( f_p \prod_{i(e)=p} (1-e^{2 \pi i \alpha_e})^{-1} \Bigr),$$
proving Theorem \ref{th:1.7}. 

To prove Theorem \ref{th:1.8}, let $f$ be an element of $K_G(\Gamma)$ of the 
form \eqref{eq:1.12} - \eqref{eq:1.13} and let 
$\phi :  V_{\Gamma} \to \RR$, $ \phi(p)=\alpha_p(\xi)$. Then
$$\chi(f) = \sum_p 
e^{2 \pi i \alpha_p}\prod_{i(e)=p} (1-e^{2 \pi i \alpha_e})^{-1}.$$
Identify $T^{\CC}$ with $\CC - 0 $ and let $C$ be a circle in the complex 
plane with radius greater than one oriented in a counter-clock-wise sense. 
Then
\begin{equation}
\label{eq:6.7}
\frac{1}{2 \pi i} \int_C \chi(f)(gz) \frac{dz}{z} = 
\sum_p \frac{1}{2 \pi i} \int_C \Bigl( 
\frac{e^{2 \pi i \alpha_p}}{\prod_{i(e)=p}(1-e^{2 \pi i \alpha_e})} \Bigr)
(gz) \frac{dz}{z}.
\end{equation}

The right hand side of this identity is easy to evaluate: 
By Proposition \ref{prop:4.5}, the summands with $\alpha_p(\xi) < 0$ 
are zero, and the summands with $\alpha_p(\xi) > 0$ are equal to
$$ Res_T \Bigl( f_p\prod_{i(e)=p} (1-e^{2 \pi i \alpha_e})^{-1} \Bigr),$$
so by Theorem \ref{th:1.7} the right hand side is equal to $\chi_{red}(f)$.
As for the the left hand side, by Theorem \ref{th:1.1}, $\chi(f)$ 
is in $R(G)$; so it is a finite sum of the form 
$$\sum m_k e^{2 \pi i \alpha_k}$$
with $m_k \in \ZZ$ and $\alpha_k \in \ZZ_G^*$; and the integral over $C$ 
of the $k$-th term is zero except when $e^{2 \pi i \alpha_k}$ doesn't 
depend on $z$, in which case the integral is just 
$2 \pi i e^{2 \pi i \alpha_k}$. Hence the left hand side is equal to 
$$\sum_{\alpha_k(\xi) =0} m_k e^{2 \pi i \alpha_k},$$
which is the character of the representation, $Q(\Gamma)^T$.

\section{GKM manifolds}
\label{sec:gkm}

Let $(M, \omega)$ be a compact $2d$-dimensional symplectic manifold 
and $\tau: G \times M \to M$ a Hamiltonian action of $G$ on $M$. We will 
say that $M$ is a \emph{symplectic GKM manifold} if $M^G$ is finite and 
if, for every $p \in M^G$, the weights 
$\alpha_{i,p} \in \ZZ_G^*$, $i=1,..,d$ of the 
isotropy representation of $G$ on $T_pM$ are pair-wise linearly independent.
Let
$$M^{(1)} = \{ p \in M ; \; \dim{G_p} \geq n-1 \}.$$
This set is called the \emph{one-skeleton} of $M$; and $M$ is a GKM 
manifold if and only if $M^{(1)}$ consists of $G$-invariant imbedded 
2-spheres, each of which contains exactly two fixed points. These 
2-spheres can intersect at the fixed points; so the combinatorial structure 
of $M^{(1)}$ is that of a graph, $\Gamma$, having the fixed points of 
$\tau$ as vertices and these 2-spheres as edges. For each oriented edge, 
$e$, of $\Gamma$, let $\varrho_e$ be the isotropy representation of $G$ on 
the tangent space to this 2-sphere at the fixed point, $t(e)$; and for each
vertex, $p$, of $\Gamma$ let $\tau_p$ be the isotropy representation of $G$
on $T_pM$. It is easily checked that $\varrho$ and $\tau$ have properties
\eqref{eq:1.2} - \eqref{eq:1.4} and hence define an action of $G$ on $\Gamma$.

For GKM manifolds the cohomology groups, $H_G(\Gamma)$ and $K_G(\Gamma)$, 
turn out to be equal to cohomology groups of $M$. More explicitly, let
$H_G(M)$ be the equivariant cohomology ring of $M$ with complex coefficients
and let $K_G(M)$ be the $K$-cohomology ring of $M$. Then there are ring
homomorphisms
\begin{eqnarray}
\label{eq:77.1}
H_G(M) & \simeq & H_G(\Gamma), \qquad (\text{see \cite{GKM}}) \\
\label{eq:77.2}
K_G(M) \otimes \CC & \simeq & K_G(\Gamma) \otimes \CC , 
\qquad (\text{see \cite{KR}}).
\end{eqnarray}
With \eqref{eq:77.1} and \eqref{eq:77.2} as our point of departure, we 
will briefly describe some geometric implications of the theorems proved 
in this paper. The first of our results, Theorem \ref{th:1.1}, is a 
``combinatorial'' explanation of why the right hand side of the 
Atiyah-Bott fixed point formula makes sense, \emph{i.e.} why 
\eqref{eq:1.10} \emph{does} define a character of a virtual 
representation of $G$. Theorems \ref{th:1.2} - \ref{th:1.4} are, in the 
manifold setting, well-known results about the ``quantum'' action of
$G$ on $M$: Suppose $[\omega] \in H^2(M, \ZZ)$. Then there exists a line 
bundle, $\LL \to M$, and a connection, $\nabla$, on this bundle with 
$curv(\nabla) =\omega$; and one says that the action, $\tau$, of $G$ on $M$ 
is \emph{pre-quantizable} if it lifts to an action of $G$ on $\LL$ preserving
$\nabla$. Now equip $M$ with a $G$-invariant Riemannian metric and let
$$ \D_{\CC} : S_{\CC}^+ \to S_{\CC}^-$$
be the $spin^{\CC}$ Dirac operator. Given the connection, $\nabla$, 
one can twist this operator with operator with $\LL$ to get a Dirac operator
$$\D_{\CC}^{\LL} : S_{\CC}^+ \otimes \LL \to S_{\CC}^- \otimes \LL,$$
and the virtual vector space
\begin{equation}
\label{eq:77.3}
Q(M) = \text{kernel}(\D_{\CC}^{\LL}) - \text{cokernel}(\D_{\CC}^{\LL})
\end{equation}
is called the $spin^{\CC}$-\emph{quantization} of $M$. 
From the action of $G$ on $\LL$, one gets a representation, $\tau_Q$, 
of $G$ on this space, and its character, $\text{trace} \tau_Q$, is equal, 
by the Atiyah-Bott formula, to the formal character, $\chi(f)$, defined by 
\eqref{eq:1.10}, $f$ being the element of $K_G(\Gamma)$ corresponding to 
$[\LL]$ under the isomorphism \eqref{eq:77.2}.

For $\tau_Q$, the convexity theorem (Theorem \ref{th:1.2}) is due to 
Guillemin and Sternberg, who pointed out in \cite{GS} that it can be 
deduced from ``quantization commutes with reduction'' and the 
Atiyah-Guillemin-Sternberg convexity theorem for moment maps. (However, 
the simple proof of this theorem described in Section \ref{sec:proof1.1} 
seems to have eluded them.) As for Theorem \ref{th:1.4}, for co-adjoint 
orbits this is the celebrated Kostant Multiplicity Theorem. Our proof of it 
in Section \ref{sec:proof1.1} is modeled on Cartier's proof of Kostant's 
theorem 
in \cite{Ca} and the symplectic version of the proof described in \cite{GLS}.

Let $T$ be a circle subgroup of $G$ and let $M_c$ be the reduction of $M$ 
with respect to $T$. For $f= [\LL]$, the ``reduced'' character, $\chi_c(f)$,
in Theorem \ref{th:1.5} can be shown, by the orbifold version of Atiyah-Bott, 
to be equal to the character of the representation of $G/T$ on $Q(M_c)$. Our
residue formula for it, (formula \eqref{eq:1.20}) appears to be a new 
result even in the manifold case; however, the formula \eqref{eq:1.21}, 
which is a special case of this formula, is just the 
``quantization commutes with reduction'' theorem for circle actions. 
A good reference for the long and entangled history of ``$[Q,R]=0$'' is the
survey article \cite{Sj}. For circle actions there are several 
relatively simple proofs, among them that of 
Duistermaat-Guillemin-Meinrenken-Wu (\cite{DGMW}), Ginzburg-Guillemin-Karshon
(\cite{GGK}) and Metzler (\cite{Me}). Of these, Metzler's proof is 
probably the closest in spirit to our combinatorial proof of Theorem 
\ref{th:1.8} in Section \ref{sec:quantred}.

\end{document}